\documentclass[12pt,a4paper]{amsart}
\usepackage[left=25mm,top=30mm,right=25mm,bottom=20mm]{geometry}
\usepackage[all]{xy}
\usepackage{amsfonts}
\usepackage{amsmath}
\usepackage{amssymb}
\usepackage{amsthm}
\usepackage{cite}
\usepackage{enumerate}
\usepackage{extarrows}
\usepackage{indentfirst}
\usepackage{mathrsfs}
\usepackage{mathtools}
\usepackage{mathptmx}
\usepackage{times}
\allowdisplaybreaks
\theoremstyle{definition}
\newtheorem{defn}{Definition}[section]
\theoremstyle{plain}
\newtheorem{thm}[defn]{Theorem}
\newtheorem{lem}[defn]{Lemma}
\newtheorem{prop}[defn]{Proposition}
\newtheorem{cor}[defn]{Corollary}
\theoremstyle{definition}
\newtheorem{eg}[defn]{Example}
\theoremstyle{definition}
\newtheorem{re}[defn]{Remark}
\def\ad{{\rm ad}}
\def\all{~\forall~}
\def\Bider{{\rm Bider}}

\def\char{{\rm char}}
\def\CBider{{\rm CBider}}
\def\Cent{{\rm Cent}}
\def\Com{{\rm Com}}
\def\CCom{{\rm CCom}}
\def\End{{\rm End}}
\def\F{{\mathbb{F}}}
\def\id{{\rm id}}
\def\N{\mathbb{N}}
\def\s{{\rm s}}
\def\set{\coloneqq}
\def\sl{{\mathfrak{sl}}}
\def\SBider{{\rm SBider}}
\def\SCom{{\rm SCom}}
\def\Z{{\mathbb{Z}}}

\newcommand{\bdefn}{\begin{defn}}
\newcommand{\edefn}{\end{defn}}
\newcommand{\bthm}{\begin{thm}}
\newcommand{\ethm}{\end{thm}}
\newcommand{\blem}{\begin{lem}}
\newcommand{\elem}{\end{lem}}
\newcommand{\bprop}{\begin{prop}}
\newcommand{\eprop}{\end{prop}}
\newcommand{\bcor}{\begin{cor}}
\newcommand{\ecor}{\end{cor}}
\newcommand{\beg}{\begin{eg}}
\newcommand{\eeg}{\end{eg}}
\newcommand{\bc}{\begin{center}}
\newcommand{\ec}{\end{center}}
\newcommand{\beq}{\begin{equation}}
\newcommand{\eeq}{\end{equation}}
\newcommand{\bre}{\begin{re}}
\newcommand{\ere}{\end{re}}
\newcommand{\bpf}{\begin{proof}}
\newcommand{\epf}{\end{proof}}

\newcommand{\benu}{\begin{enumerate}}
\newcommand{\eenu}{\end{enumerate}}
\newcommand{\Beq}{\begin{equation*}}
\newcommand{\Eeq}{\end{equation*}}
\newcommand{\bspl}{\begin{split}}
\newcommand{\espl}{\end{split}}
\newcommand\relphantom[1]{\mathrel{\phantom{#1}}}

\numberwithin{equation}{section}
\begin{document}
\title[Biderivations and commuting linear maps on Hom-Lie algebras]{Biderivations and commuting linear maps\\ on Hom-Lie algebras}
\author{Bing Sun, Yao Ma and Liangyun Chen}
\address{B. Sun: School of Mathematics, Changchun Normal
 University, Changchun 130032, China}
\address{Y. Ma: School of
Mathematics and Statistics, Northeast Normal
 University, Changchun 130024, China}
\address{L. Chen: School of
Mathematics and Statistics, Northeast Normal
 University, Changchun 130024, China}

\thanks{*Corresponding author (Y. Ma): may703@nenu.edu.cn.}
\thanks{Supported by NNSF of China
(Nos. 11901057, 11801066, 11771069 and 11771410), NSF of Jilin province (No.
20170101048JC), Natural Science Foundation of Changchun Normal University. }
\begin{abstract}
The purpose of this paper is to determine skew-symmetric biderivations $\Bider_\s(L, V)$
and commuting linear maps $\Com(L, V)$ on a Hom-Lie algebra $(L,\alpha)$ having their ranges in an $(L,\alpha)$-module $(V, \rho, \beta)$, which are both closely related to $\Cent (L, V)$, the centroid of $(V, \rho, \beta)$. Specifically, under appropriate assumptions, every $\delta\in\Bider_\s(L, V)$ is of the form $\delta(x,y)=\beta^{-1}\gamma([x,y])$ for some $\gamma\in \Cent (L, V)$, and $\Com(L, V)$ coincides with $\Cent (L, V)$. Besides, we give the algorithm for describing $\Bider_\s(L, V)$ and $\Com(L, V)$ respectively, and provide several examples.
\bigskip

\noindent{\bf Keywords:} Biderivation, commuting linear map, centroid, Hom-Lie algebra.\\
\noindent{\bf MSC(2020):} 17B40, 17B61, 16R60.
\end{abstract}
\maketitle
 \section{Introduction}

The notion of biderivations appeared in different areas. Maksa used biderivations to study real Hilbert space \cite{Maksa}. Vukman investigated symmetric biderivations in the prime and semiprime rings \cite{Vukman}. The well-known result that every biderivation on a noncommutative prime ring $A$ is of the form $\lambda[x,y]$ for some $\lambda$ belonging to the extended centroid of $A$, was discovered independently by Bresar et al \cite{Bresar&Martindale&Miers}, Skosyrskii \cite{Skosyrskii}, and Farkas and Letzter \cite{Farkas&Letzter}, where biderivations were connected with noncommutative Jordan algebras by Skosyrskii and with Poisson algebras by Farkas and Letzter, respectively. Besides their wide applications, biderivations are interesting in their own right and have been introduced to Lie algebras \cite{Wang&Yu&Chen}, which were studied by many authors recently. In particular, biderivations are closely related to the theory of commuting linear maps that has a long and rich history, and we refer to the survey \cite{Bresar} for the development of commuting maps and their applications. It is worth mentioning that Bre\v{s}ar and Zhao considered a general but simple approach for describing biderivations and commuting linear maps on a Lie algebra $L$ that having their ranges in an $L$-module \cite{Bresar&Zhao}, which covered most of the results in \cite{Chen, Liu&Guo&Zhao, Han&Wang&Xia, Tang, Wang&Yu}, and inspires us to generalize their method to Hom-Lie algebras.

Hom-Lie algebras are a generalization of Lie algebras, the notion of which was initially introduced in \cite{Hartwig&Larsson&Silvestrov}, motivated by the study of quantum deformations or $q$-deformations of the Witt and the Virasoro algebras via twisted derivations, after several investigations of various quantum deformations of Lie algebras\cite{Aizawa&Sato, Chaichian&Kulish, Curtright&Zachos, Hu, Kassel}. Due to the close relation to discrete and deformed vector fields and differential calculus, this kind of algebraic structures have been studied extensively during the last decade, which we mention here just a few: representations, (co)homology and deformations of Hom-Lie algebras \cite{Agrebaoui&Benali&Makhlouf, Alvarez&Cartes, Makhlouf&Silvestrov, Sheng, Yau}; structure theory of simple and semisimple Hom-Lie algebras \cite{Chen&Han, Jin&Li, Xie&Liu}; Hom-Lie structures on Kac-Moody algebras \cite{Makhlouf&Zusmanovich}; extensions of Hom-Lie algebras \cite{Casas&GarciaMartinez, Larsson&Silvestrov}; geometric generalization of Hom-Lie algebras \cite{Laurent-Gengoux&Teles, Peyghan&Nourmohammadifar}; integration of Hom-Lie algebras \cite{Laurent-Gengoux&Makhlouf&Teles}.

Hence it would be natural to generalize known theories from Lie algebras to Hom-Lie algebras, and we are interested in determining biderivations and commuting linear maps on a Hom-Lie algebra follows from \cite{Bresar&Zhao}. The paper is organized as follows.

In Section 2, after recalling some preliminaries on Hom-Lie algebras, we introduce the notion of Hom-Lie algebra module homomorphisms, and give the Schur's lemma for adjoint Hom-Lie algebra modules.

In Section 3, we define biderivations on a Hom-Lie algebra $(L,\alpha)$ having their ranges in an $(L,\alpha)$-module $(V, \rho, \beta)$ as well as the centroid of $(V, \rho, \beta)$, and derive that
skew-symmetric biderivations arise from the centroid provided that
$(L,\alpha)$ is perfect, $\alpha$ is surjective, $\beta$ is invertible and $Z_V(L)=\{0\}$ (see Theorem \ref{BiderofperfectHomLiealg}). In particular, every skew-symmetric biderivation on a simple $(L,\alpha)$ with $\alpha$ invertible having its range in the $\alpha^k$-adjoint $(L,\alpha)$-module is of the form $\lambda\alpha^k([x,y])$ (see Theorem \ref{BiderofsimpleHomLiealg}). We also give an algorithm to find all skew-symmetric biderivations and apply it to several examples.

In Section 4, we give the definition of commuting linear maps from a Hom-Lie algebra $(L,\alpha)$ to an $(L,\alpha)$-module $(V, \rho, \beta)$, which coincides with the centroid of $(V, \rho, \beta)$ if $\alpha$ is surjective, $\beta$ is invertible and $Z_V(L')=\{0\}$ (see Theorem \ref{ComofHomLiealg}). An
algorithm to describe all commuting linear maps on $(L,\alpha)$ is also provided.

Throughout this paper, all the vector spaces are over a fixed field $\F$ such that $\char(\F)\neq 2$ unless otherwise stated.

\section{The Schur's
lemma for adjoint modules over a Hom-Lie algebra}

We begin with some definitions concerning Hom-Lie algebras.
\bdefn\cite{Sheng}
A \textit{Hom-Lie algebra} $(L,\alpha)$ is an $\F$-vector space $L$ endowed with a bilinear map $[-,-]: L\times L \rightarrow L$ and a linear homomorphism $\alpha: L\rightarrow L$ satisfying that for any $x, y, z \in L$,
\begin{gather*}
[x,y]=-[y,x], ~~\text{(skew-symmetry)}\\
[\alpha(x),[y,z]]+[\alpha(y),[z,x]]+[\alpha(z),[x,y]]=0. ~~\text{(Hom-Jacobi identity)}
\end{gather*}
In particular, a Hom-Lie algebra $(L,\alpha)$ is called \textit{multiplicative}, if $\alpha$ is an algebra homomorphism, i.e., $\alpha([x,y])=[\alpha(x),\alpha(y)]$, for all $x, y\in L$.
\edefn

For the structural study of $(L,\alpha)$ in this paper, we need the stability properties that the multiplicative condition of $\alpha$ offers, so in the sequel, \textbf{we will say Hom-Lie algebras when referring to multiplicative Hom-Lie algebras}.

Moreover, a subspace $I$ of $(L,\alpha)$ is called an \textit{ideal}, if $[IL]\subseteq I$ and $\alpha(I)\subseteq I$. If $(L,\alpha)$ has no proper ideals and is not abelian, then we call it \textit{simple}. $(L,\alpha)$ is called \textit{perfect} if $L'\set [L,L]=L$. The \textit{center} of a Hom-Lie algebra $(L,\alpha)$, denoted by $Z(L)$, is the set $\{z\in L\ |\ [z,L]=0\}$. For more details of these definitions, we refer to \cite{Casas&GarciaMartinez, Chen&Han, Makhlouf&Zusmanovich}.

\bdefn\cite{Sheng,Yau}
For a Hom-Lie algebra $(L,\alpha)$, a triple $(V, \rho, \beta)$ consisting of a vector space $V$, a linear map
$\rho: L\rightarrow \End(V)$ and $\beta\in \End(V)$ is said to be a \textit{representation} of $(L,\alpha)$ or an $(L,\alpha)$-\textit{module}, if for all $x,y\in L$, the following equalities are satisfied
\begin{gather}
\beta\circ \rho(x)=\rho(\alpha(x))\circ\beta,\label{bp=pab}\\
\rho([x,y])\circ\beta=\rho(\alpha(x))\circ\rho(y)-\rho(\alpha(y))\circ\rho(x).\label{pb=pap-pap}
\end{gather}
For brevity of notation, we usually put $xv=\rho(x)(v), \all x\in L, v\in V$, just like the case in Lie algebras. A subspace $W$ of $V$ is called a \textit{submodule} of $(V, \rho, \beta)$ if $W$ is both $\beta$-invariant and $L$-invariant, i.e., $\beta(W)\subset W$ and $xW\subset W, \all x\in L$.
\edefn

\beg
(1) For any integer $k$, $(L,\ad_k,\alpha)$ with $\ad_k(x)(y)\set [\alpha^k(x),y]$ is an $(L,\alpha)$-module by \cite[Lemma 6.2]{Sheng}, called the $\alpha^k$-adjoint $(L,\alpha)$-module. Note that $\alpha$ need to be invertible in case $k$ is negative.

(2) If $(V, \rho, \beta)$ is an $(L,\alpha)$-module, then $(V, \rho_k\set \rho \alpha^k, \beta)$ is also an $(L,\alpha)$-module, which follows from the equalities that
\begin{align*}
\beta \rho_k(x)(v) &=\beta\rho(\alpha^k(x))(v) \overset{(\ref{bp=pab})}{=}\rho(\alpha^{k+1}(x))\beta(v) =\rho_k(\alpha(x))\beta(v),\\
\rho_k([x,y])\beta(v) &=\rho\alpha^k([x,y])\beta(v) =\rho([\alpha^k(x),\alpha^k(y)])\beta(v) \\
\text{(by (\ref{pb=pap-pap}))}&=\rho(\alpha^{k+1}(x))\rho(\alpha^k(y))(v)-\rho(\alpha^{k+1}(y))\rho(\alpha^k(x))(v)\\
&=\rho_k(\alpha(x))\rho_k(y)(v)-\rho_k(\alpha(y))\rho_k(x)(v).
\end{align*}
\eeg

\bre\label{re1}
If $I$ is a submodule of $(L,\ad_k,\alpha)$ and $\alpha^k(L)=L$, then $I$ is an ideal of $L$. In fact,
\Beq
\alpha(I)\subset I ~~{\rm and}~~ [L,I]=[\alpha^k(L),I]=\ad_k(L)(I)\subset I.
\Eeq
\ere

Next we introduce the notion of homomorphisms between representations of a Hom-Lie algebra $(L,\alpha)$ and show that the Schur's Lemma holds for finite-dimensional adjoint $(L,\alpha)$-modules.
\bdefn
Let $(V_1,\rho_1,\beta_1)$ and $(V_2,\rho_2,\beta_2)$ be two $(L,\alpha)$-modules. A linear map $f: V_1\rightarrow V_2$ is said to be a \textit{homomorphism} of $(L,\alpha)$-modules, if
$\beta_2 \circ f=f \circ \beta_1$ and $f\circ \rho_1(x)=\rho_2(\alpha(x))\circ f, \all x\in L$, or, in terms of two commutative diagrams,
\Beq
\xymatrix{
 V_1 \ar[d]_{f} \ar[r]^{\beta_1} & V_1 \ar[d]^{f} \\
 V_2 \ar[r]_{\beta_2} & V_2,
 }
\qquad
\xymatrix{
 V_1 \ar[d]_{f} \ar[r]^{\rho_1(x)} & V_1 \ar[d]^{f} \\
 V_2 \ar[r]_{\rho_2\alpha(x)} & V_2.
 }
\Eeq

\edefn

Note that the identity map is no longer a natural homomorphism on a Hom-Lie algebra module in general, which is actually twisted by $\alpha$ in some way.

\beg\label{homoofadmod}
For any integers $k, s$, it follows that
$\alpha^{s+1}:(L,\ad_k,\alpha)\rightarrow (L,\ad_{k+s},\alpha)$ is a homomorphism of $(L,\alpha)$-modules, since $\alpha\circ \alpha^{s+1}=\alpha^{s+1} \circ \alpha$ and
\begin{gather*}
\alpha^{s+1}{\ad}_k(x)(y) =\alpha^{s+1}[\alpha^{k}(x),y] =[\alpha^{k+s+1}(x),\alpha^{s+1}(y)]={\ad}_{k+s}(\alpha(x))\alpha^{s+1}(y).
\end{gather*}
In particular, $\alpha^0=\id_L$ is a homomorphism from $(L,\ad_0,\alpha)$ to $(L,\ad_{-1},\alpha)$ in the case that $\alpha$ is invertible.
\eeg

By direct calculations, we have the following property.
\bprop
Let $(V_1,\rho_1,\beta_1)$ and $(V_2,\rho_2,\beta_2)$ be two $(L,\alpha)$-modules and $f:V_1\rightarrow V_2$ a module homomorphism. Then ${\ker}f$ is a submodule of $(V_1,\rho_1,\beta_1)$.
\eprop

\bthm[Schur's lemma] \label{SchurLemma}
Let $(L,\alpha)$ be a Hom-Lie algebra with $\alpha$ invertible, which is called \textit{regular} in \cite{Sheng}. Let the base field $\F$ be algebraically closed.
If $(L,\alpha)$ is simple and finite-dimensional, and $f:(L,\ad_k,\alpha)\rightarrow (L,\ad_{k+s},\alpha)$ is a homomorphism of $(L,\alpha)$-modules, then
\Beq
f=\lambda \alpha^{s+1}, ~~\text{for some~} \lambda\in \F.
\Eeq
\ethm

\bpf
We know that $f\alpha^{-s-1}:L\rightarrow L$ is a linear map. Then use the facts that $(L,\alpha)$ is finite-dimensional and $\F$ is algebraically closed to find an eigenvector $0\neq x_0\in L$ of $f\alpha^{-s-1}$ for some eigenvalue, say $\lambda\in \F$, i.e., $f\alpha^{-s-1}(x_0)=\lambda x_0$. Since $f\circ\alpha=\alpha \circ f$, $f(x_0)=\lambda \alpha^{s+1}(x_0)$.

By Example \ref{homoofadmod}, $\lambda \alpha^{s+1}:(L,\ad_k,\alpha)\rightarrow (L,\ad_{k+s},\alpha)$ is an $(L,\alpha)$-module homomorphism, which implies that $f-\lambda\alpha^{s+1}$ is also an $(L,\alpha)$-module homomorphism. Then $\ker(f-\lambda\alpha^{s+1})$ is a non-zero submodule of $(L,\ad_k,\alpha)$, so is a non-zero ideal of $(L,\alpha)$, by Remark \ref{re1}. Hence $\ker(f-\lambda\alpha^{s+1})=L$ by the simplicity of $(L,\alpha)$, that is, $f=\lambda \alpha^{s+1}$.
\epf

\section{Biderivations of Hom-Lie algebras}

Let $(V, \rho, \beta)$ be an $(L,\alpha)$-module. A bilinear map $\delta: L\times L\rightarrow V$ is said to be a \textit{biderivation} if
\begin{gather}
\beta\delta(x, y)=\delta(\alpha(x), \alpha(y)),\label{bd=da}\\
\delta(\alpha(z),[x,y])=\alpha(x)\delta(z,y)-\alpha(y)\delta(z,x),\label{fixleftder}\\
\delta([x,y],\alpha(z))=\alpha(x)\delta(y,z)-\alpha(y)\delta(x,z),\label{fixrightder}
\end{gather}
for all $x, y, z\in L$. We say that a biderivation $\delta$ is \textit{skew-symmetric}, if
\beq\label{skewsymm}
\delta(x,y)=-\delta(y,x).
\eeq
If this is the case, Eqs. (\ref{fixleftder}) and (\ref{fixrightder}) coincide. Hence, a bilinear map $\delta: L\times L\rightarrow V$ is a skew-symmetric biderivation provided that $\delta$ satisfies (\ref{skewsymm}), (\ref{bd=da}), and any one of (\ref{fixleftder}) or (\ref{fixrightder}).

Denote by $\Bider(L, V)$ (resp. $\Bider_\s(L, V)$) the set of all biderivations (resp. skew-symmetric biderivations) $\delta: L\times L\rightarrow V$. In particular, we write $\Bider_\s(L, \ad_k)$ as the set of all skew-symmetric biderivations on $(L, \alpha)$ for the adjoint module $(L, \ad_k, \alpha)$.

Using the Hom-Jacobi identity, it is easy to see that a skew-symmetric biderivation $\delta: L\times L\rightarrow V$ is a generalization of the Hom-Lie bracket $[-,-]: L\times L\rightarrow L$, where the last $L$ is viewed as $(L, \ad_0, \alpha)$.

Now we would like to explain the terminology ``biderivation''. Note that the notion of derivations on a Hom-Lie algebra $(L,\alpha)$ was introduced in \cite{Sheng}, where a linear map $D: L\rightarrow L$ is called an $\alpha^k$-\textit{derivation} of $(L,\alpha)$, if $D \circ\alpha=\alpha\circ D$ and
\Beq
D[x,y]=[\alpha^k(x), D(y)]-[\alpha^k(y),D(x)], ~\all x, y\in L.
\Eeq
Under the above definition, for any nonnegative integer $k$ and any $z\in L$ satisfying $\alpha(z)=z$, a linear map $D_k(z): L\rightarrow L$ defined by
\Beq
D_k(z)(x)=[z, \alpha^k(x)], ~\all x\in L
\Eeq
becomes an $\alpha^{k+1}$-derivation, which implies that $[z,-]$ and $[-,z]$ are $\alpha$-derivations, then it is reasonable to call $[-,-]$ a biderivation. Finally, its generalization $\delta$ can be called by the same name.

We've already seen that $[-,-]\in\Bider_\s(L, \ad_0)$. To provide more examples of biderivations, we define the \textit{centroid} of $(V, \rho, \beta)$ as
\Beq
\Cent (L, V)=\{\gamma:L\rightarrow V \mid \gamma([x,y])=\alpha(x)\gamma(y)~~{\rm and}~~\beta\circ \gamma=\gamma\circ \alpha\}.
\Eeq
Clearly, $\Cent (L, V)$ is just the space of $(L,\alpha)$-module homomorphisms from $(L, \ad_0, \alpha)$ to $(V, \rho, \beta)$. In particular, we write $\Cent (L, \ad_k)$ as the centroid of the adjoint module $(L, \ad_k, \alpha)$.
\beg
Let $(V, \rho, \beta)$ be a representation of a Hom-Lie algebra
$(L,\alpha)$ with $\beta$ invertible. If $\gamma\in \Cent (L, V)$, then
\Beq
\delta(x,y)=\beta^{-1}\gamma([x,y]),~\all x, y\in L,
\Eeq
is a skew-symmetric biderivation.

In fact, $\delta$ is clearly skew-symmetric; (\ref{bd=da}) holds since
\Beq
\delta(\alpha(x), \alpha(y))=\beta^{-1}\gamma([\alpha(x), \alpha(y)])=\beta^{-1}\gamma\alpha([x,y])=\gamma([x,y])=\beta\delta(x,y), ~\all x, y\in L;
\Eeq
(\ref{fixrightder}) holds since
\begin{align*}
\delta([x,y],\alpha(z))&=\beta^{-1}\gamma([[x,y],\alpha(z)])\\ \text{(by Hom-Jacobi identity)}&=\beta^{-1}\gamma([\alpha(x),[y,z]])-\beta^{-1}\gamma([\alpha(y),[x,z]])\\
&=\beta^{-1}(\alpha^2(x)\gamma([y,z])-\beta^{-1}\alpha^2(y)\gamma([x,z])\\
\text{(by (\ref{bp=pab}))}&=\alpha(x)\beta^{-1}\gamma([y,z])-\alpha(y)\beta^{-1}\gamma([x,z])\\
&=\alpha(x)\delta(y,z)-\alpha(y)\delta(x,z).
\end{align*}
\eeg

From the example above we could see that for an $(L,\alpha)$-module $(V, \rho, \beta)$ such that $\beta$ is invertible, each $\gamma\in \Cent (L, V)$ induces a skew-symmetric biderivation $\delta(x,y)=\beta^{-1}\gamma([x,y])$. Later, we will show that, under appropriate assumptions, all biderivations are of this form, i.e., Theorem \ref{BiderofperfectHomLiealg}, to prove which we need the following technical lemmas.

In Lemmas \ref{xyzw=0}, \ref{cyclicsumxyz}, \ref{d(u,[xy])-ud(x,y)} and \ref{d(z,[xy])=zd(x,y)}, $(V, \rho, \beta)$ is a representation of a Hom-Lie algebra
$(L,\alpha)$ and $\delta\in\Bider_\s(L, V)$. Moreover, we denote by
\Beq
Z_V(S)=\{v\in V\ |\ sv=0,~\all s\in S\}
\Eeq
the subspace of $V$ that is killed by any subset $S$ of $L$. Note that for an adjoint module $(L,\ad_k,\alpha)$ with $\alpha$ surjective,
$Z_{L}(L)$ is an ideal of $(L,\alpha)$ that coincides with $Z(L)$, the center of $(L,\alpha)$.

\blem\label{xyzw=0}
$\beta([x,y]\delta(z,w)+[z,w]\delta(x,y))=0, \all x,y,z,w\in L.$
\elem

\bpf
Consider the following quaternary linear map
\begin{align*}
\varphi:L\times L\times L \times L&\rightarrow V,\\
(x,y,z,w)&\mapsto [x,y]\delta(z,w)+[z,w]\delta(x,y).
\end{align*}
Then it suffices to prove $\beta\varphi(x,y,z,w)=0, \all x,y,z,w\in L$.

Clearly, $\varphi$ satisfies $\varphi(x,y,z,w)=\varphi(z,w,x,y)$ and
\beq\label{12}
\varphi(x,y,z,w)=-\varphi(y,x,z,w)
\eeq
by the definition.

We'll prove $\beta\varphi(x,y,z,w)=\beta\varphi(z,y,x,w)$, by computing $\delta([x,y],[z,w])$ in two ways. For any $x,y,z,w\in L$, we have
\begin{align*}
&\delta(\alpha([x,y]),\alpha([z,w]))\\
=&\delta([\alpha(x), \alpha(y)],\alpha([z,w]))\\
=&\alpha^2(x)\delta(\alpha(y),[z,w])-\alpha^2(y)\delta(\alpha(x),[z,w])\\
=&\alpha^2(x)\big(\alpha(z)\delta(y,w)-\alpha(w)\delta(y,z)\big)-\alpha^2(y)\big(\alpha(z)\delta(x,w)-\alpha(w)\delta(x,z)\big)\\
=&\alpha^2(x)\alpha(z)\delta(y,w)-\alpha^2(x)\alpha(w)\delta(y,z)-\alpha^2(y)\alpha(z)\delta(x,w)+\alpha^2(y)\alpha(w)\delta(x,z).
\end{align*}
On the other hand,
\begin{align*}
&\delta(\alpha([x,y]),\alpha([z,w]))\\
=&\delta(\alpha([x,y]),[\alpha(z), \alpha(w)])\\
=&\alpha^2(z)\delta([x,y],\alpha(w))-\alpha^2(w)\delta([x,y],\alpha(z))\\
=&\alpha^2(z)\big(\alpha(x)\delta(y,w)-\alpha(y)\delta(x,w)\big)-\alpha^2(w)\big(\alpha(x)\delta(y,z)-\alpha(y)\delta(x,z)\big)\\
=&\alpha^2(z)\alpha(x)\delta(y,w)-\alpha^2(z)\alpha(y)\delta(x,w)-\alpha^2(w)\alpha(x)\delta(y,z)+\alpha^2(w)\alpha(y)\delta(x,z).
\end{align*}
Comparing above relations and using (\ref{pb=pap-pap}), we have
\Beq
[\alpha(x),\alpha(z)]\beta \delta (y,w)+[\alpha(y),\alpha(w)]\beta \delta (x,z)=[\alpha(y),\alpha(z)]\beta \delta (x,w)
+[\alpha(x),\alpha(w)]\beta \delta (y,z).
\Eeq
By (\ref{bp=pab}), it follows that $\beta\varphi(x,z,y,w)=\beta\varphi(y,z,x,w)$, which implies
\beq \label{13}
\beta\varphi(x,y,z,w)=\beta\varphi(z,y,x,w), ~\all x,y,z,w\in L.
\eeq
Then
\beq \label{23}
\beta\varphi(x,y,z,w) \overset{(\ref{13})}{=}\beta\varphi(z,y,x,w) \overset{(\ref{12})}{=}-\beta\varphi(y,z,x,w) \overset{(\ref{13})}{=}-\beta\varphi(x,z,y,w).
\eeq
It follows that
\Beq
\beta\varphi(x,y,z,w) \overset{(\ref{12})}{=}-\beta\varphi(y,x,z,w) \overset{(\ref{23})}{=}\beta\varphi(y,z,x,w) \overset{(\ref{12})}{=}-\beta\varphi(z,y,x,w) \overset{(\ref{13})}{=}-\beta\varphi(x,y,z,w).
\Eeq
The proof is complete.
\epf

\blem\label{cyclicsumxyz}
$\circlearrowleft_{x,y,z}\delta([x,y],\alpha(z))=2(\alpha(z)\delta(x,y)-\delta(\alpha(z),[x,y])), \all x, y, z\in L$, where $\circlearrowleft_{x,y,z}$ denotes the sum over cyclic permutations of $x, y, z$.
\elem
\bpf
Note that
\begin{align*}
&\circlearrowleft_{x,y,z}\delta([x,y],\alpha(z))\\
=&\delta([x,y],\alpha(z))+\delta([y,z],\alpha(x))+\delta([z,x],\alpha(y))\\
=&\alpha(x)\delta(y,z)-\alpha(y)\delta(x,z)+\alpha(y)\delta(z,x)-\alpha(z)\delta(y,x)+\alpha(z)\delta(x,y)-\alpha(x)\delta(z,y)\\
=&2\big(\alpha(x)\delta(y,z)-\alpha(y)\delta(x,z)+\alpha(z)\delta(x,y)\big)\\
=&2\big(\alpha(z)\delta(x,y)-\delta(\alpha(z),[x,y])\big).
\end{align*}
Hence the lemma follows.
\epf

\blem\label{d(u,[xy])-ud(x,y)}
If $\beta$ is invertible, then
\Beq
\delta(\alpha(u),[x,y])-\alpha(u)\delta(x,y)\in Z_V(L'), ~\all x,y,u\in L.
\Eeq
\elem

\bpf
By Lemma \ref{xyzw=0} and (\ref{bp=pab}), we have
\Beq [\alpha(x),\alpha(y)]\beta\delta(z,w)+[\alpha(z),\alpha(w)]\beta\delta(x,y)=0.
\Eeq Then
\begin{align*}
0=&\circlearrowleft_{x,u,y}\big([[\alpha(x),\alpha(u)],\alpha^2(y)]\beta\delta(z,w)+[\alpha(z),\alpha(w)]\beta\delta([x,u],\alpha(y))\big)\\
=&\circlearrowleft_{x,u,y}\big[[\alpha(x),\alpha(u)],\alpha^2(y)]\beta\delta(z,w)+\circlearrowleft_{x,u,y}[\alpha(z),\alpha(w)]\beta\delta([x,u],\alpha(y))\\
=&-[\alpha(z),\alpha(w)]\beta\big(\circlearrowleft_{x,y,u}\delta([x,y],\alpha(u))\big)\\
=&2\beta\big([z,w](\delta(\alpha(u),[x,y])-\alpha(u)\delta(x,y)\big),
\end{align*}
where the last equality uses Lemma \ref{cyclicsumxyz} and (\ref{bp=pab}). Hence the lemma follows from the fact that $\beta$ is invertible and $\char\F\neq2$.
\epf

\blem\label{d(z,[xy])=zd(x,y)}
If $\beta$ is invertible and $(L,\alpha)$ is a perfect Hom-Lie algebra with $\alpha$ surjective, then $\delta(z,[x,y])=z\delta(x,y), \all x,y,z\in L.$
\elem

\bpf
Since $(L,\alpha)$ is a perfect Hom-Lie algebra with $\alpha$ surjective, for any $x,y,z\in L$, there exist $u_i,v_i(1\leq i\leq m),s,t\in L$ such that $z=\sum_{i=1}^{m}[\alpha(u_i),\alpha(v_i)]$ and $[x,y]=[\alpha(s),\alpha(t)]$. Then
\begin{align*}
&\delta(z,[x,y])-z\delta(x,y)\\
=&\sum_{i=1}^{m}\big(\delta([\alpha(u_i),\alpha(v_i)],[\alpha(s),\alpha(t)]) -[\alpha(u_i),\alpha(v_i)]\delta(\alpha(s),\alpha(t))\big)\\
=&\sum_{i=1}^{m}\big(\delta([\alpha(u_i),\alpha(v_i)],\alpha([s,t])) -[\alpha(u_i),\alpha(v_i)]\beta\delta(s,t)\big)\\
=&\sum_{i=1}^{m}\big(\alpha^2(u_i)\delta(\alpha(v_i),[s,t]) -\alpha^2(v_i)\delta(\alpha(u_i),[s,t]) -\alpha^2(u_i)\alpha(v_i)\delta(s,t)+\alpha^2(v_i)\alpha(u_i)\delta(s,t)\big)\\
=&\sum_{i=1}^{m}\big(\alpha^2(u_i)(\delta(\alpha(v_i),[s,t]) -\alpha(v_i)\delta(s,t)) -\alpha^2(v_i)(\delta(\alpha(u_i),[s,t]) -\alpha(u_i)\delta(s,t))\big).
\end{align*}
By Lemma \ref{d(u,[xy])-ud(x,y)}, we obtain that the last equation above is equal to zero.
\epf

\bthm\label{BiderofperfectHomLiealg}
Let $(L,\alpha)$ be a perfect Hom-Lie algebra with $\alpha$ surjective and $(V, \rho, \beta)$ an $(L,\alpha)$-module with $\beta$ invertible such that $Z_V(L)={0}$. Then for all $\delta\in\Bider_\s(L, V)$, there exists $\gamma\in \Cent (L, V)$ such that $\delta(x,y)=\beta^{-1}\gamma([x,y]), \all x,y\in L$.
\ethm
\bpf
Since $L=L'$, for any $z\in L$, we have $z=\sum_{i=1}^{m}[z_i',z_i'']$ for some $z_i',z_i''\in L.$
Define $\gamma: L\rightarrow V$ by
\Beq \gamma\bigg(\sum_{i=1}^{m}[z_i',z_i'']\bigg)=\sum_{i=1}^{m}\beta\delta(z_i',z_i'').
\Eeq Then it suffices to show that $\gamma$ is well-defined as well as $\gamma([x,y])=\alpha(x)\gamma(y), \all x,y\in L.$ Indeed,

If $\sum_{i=1}^{m}[z_i',z_i'']=0$, then for any $u\in L$,
\Beq 0 =\delta\Big(u,\sum_{i=1}^{m}[z_i',z_i'']\Big) =\sum_{i=1}^{m}\delta(u,[z_i',z_i'']) =u\sum_{i=1}^{m}\delta(z_i',z_i''),
\Eeq
and hence $\sum_{i=1}^{m}\delta(z_i',z_i'')=0$ follows from $Z_V(L)={0}$. So $\sum_{i=1}^{m}\beta\delta(z_i',z_i'')=0$, which implies that $\gamma$ is well-defined. Furthermore, suppose $y_i',y_i''\in L$ such that $y=\sum_{i=1}^{k}[y_i',y_i'']$, then
\begin{align*}
\gamma([x,y])=
&\beta\delta(x,y) =\delta(\alpha(x),\alpha(y)) =\delta\Big(\alpha(x),\sum_{i=1}^{k}[\alpha(y_i'),\alpha(y_i'')]\Big) =\sum_{i=1}^{k}\alpha(x)\delta(\alpha(y_i'),\alpha(y_i''))\\
=&\sum_{i=1}^{k}\alpha(x)\beta\delta(y_i',y_i'') =\sum_{i=1}^{k}\alpha(x)\gamma([y_i',y_i'']) =\alpha(x)\gamma(y).
\end{align*}
Therefore, $\gamma\in \Cent (L, V)$.
\epf

Now we turn to $\Bider_\s(L, \ad_k)$, which is easy to determine in the case that $(L,\alpha)$ is centerless and perfect, or in particular, $(L,\alpha)$ is simple.

\bthm \label{BiderofsimpleHomLiealg}
If $(L,\alpha)$ is a centerless and perfect Hom-Lie algebra with $\alpha$ invertible, then every skew-symmetric biderivation $\delta:L\times L\rightarrow (L,\ad_k,\alpha)$ is of the form
\Beq
\delta(x,y)=\alpha^{-1}\gamma([x,y]),
\Eeq
where $\gamma\in \Cent (L, \ad_k)$. Moreover, if $(L,\alpha)$ is simple and finite-dimensional, then
\Beq
\delta(x,y)=\lambda\alpha^{k}([x,y]), ~~\text{for some~} \lambda\in \F.
\Eeq
\ethm
\bpf
The first assertion is an immediate corollary to Theorem \ref{BiderofperfectHomLiealg} when $(L,\alpha)$ is centerless and perfect. Now let $(L,\alpha)$ be simple and finite-dimensional.

If $\delta=0$, then take $\lambda=0$ and we are done.

Suppose $\delta\neq 0$. Using Theorem \ref{BiderofperfectHomLiealg}, there exists $\gamma\in \Cent (L, \ad_k)$ such that $\delta(x,y)=\alpha^{-1}\gamma([x,y])$. Note that $\gamma:(L,\ad_0,\alpha)\rightarrow (L,\ad_k,\alpha)$ is a module homomorphism over a finite-dimensional simple Hom-Lie algebra $(L,\alpha)$ together with the fact that $\alpha$ is invertible. By Theorem \ref{SchurLemma}, $\gamma=\lambda \alpha^{k+1}$ and $\delta(x,y)=\lambda\alpha^{k}([x,y]).$
\epf

In general, on a Hom-Lie algebra $(L,\alpha)$ such that $\alpha$ is invertible, we would like to give an algorithm to find all skew-symmetric biderivations in $\Bider_\s(L, \ad_k)$. To this end, we need the notions of central and special biderivations.

\bdefn \label{CBider&SBider}
Let $(L,\alpha)$ be a Hom-Lie algebra and $\delta\in \Bider_\s(L, \ad_k)$. $\delta$ is called \textit{central} if
\Beq
\delta(L,L)\subset Z(L)=Z_L(L);
\Eeq
$\delta$ is called \textit{special} if
\Beq
\delta(L,L)\subset Z_L(L')~~\text{and}~~ \delta(L',L')=0.
\Eeq
Denote by $\CBider_\s(L, \ad_k)$ and $\SBider_\s(L, \ad_k)$ the sets consisting of all central and special skew-symmetric biderivations in $\Bider_\s(L, \ad_k)$, respectively.
\edefn

\bre \label{CBider}
Clearly, each $\delta\in\CBider_\s(L, \ad_k)$ satisfies $\delta(L,L')=0$. Then it follows that $\CBider_\s(L, \ad_k)\subseteq \SBider_\s(L, \ad_k)$, and any skew-symmetric bilinear map $\delta$ satisfying
\begin{gather*}
\beta\delta(x, y)=\delta(\alpha(x), \alpha(y)), ~\all x, y\in L,\\
\delta(L, L')=0,\\
\delta(L, L)\subseteq Z(L)
\end{gather*}
automatically belongs to $\CBider_\s(L, \ad_k)$.
\ere

\beg
Every Hom-Lie algebra $(L,\alpha)$ with nontrivial center and such that the codimension of $L'$ in $L$ is no less than $2$ has nonzero special biderivations. Indeed, by the codimension assumption, there exists a nonzero skew-symmetric bilinear form $\omega : L \times L \rightarrow \F$ such that $\omega(L,L')={0}$, and taking any nonzero $z_0\in Z(L)$, we have that $\delta(x,y)\set \omega(x,y)z_0$ is a nonzero central biderivation as well as a nonzero special biderivation.
\eeg

Suppose that $(L,\alpha)$ is a Hom-Lie algebra with $\alpha$ surjective and $\delta\in\Bider_\s(L,\ad_k)$. The motivation for our algorithm is the following observation.

Note that for any $x,y\in L$, $z\in Z(L)$ we have
\Beq
0=\delta([z,x],\alpha(y)) =[\alpha^{k+1}(z),\delta(x,y)] -[\alpha^{k+1}(x),\delta(z,y)] =-[\alpha^{k+1}(x),\delta(z,y)],
\Eeq
which implies $\delta(Z(L),L)\subset Z(L)$. Since $Z(L)$ is an ideal of $(L,\alpha)$, we consider the quotient Hom-Lie algebra $(\bar{L},\bar{\alpha})$ and its quotient adjoint module $(\bar{L},\overline{\ad}_k,\bar{\alpha})$, where $\bar{L}\set L/Z(L)$, $\bar{\alpha}(\bar{x})\set \overline{\alpha(x)},$ and $\overline{\ad}_k(\bar{x})(\bar{y})\set \overline{[\alpha^k(x), y]}$, for all $x, y\in L$. Then we can define a skew-symmetric biderivation $\bar{\delta}\in\Bider_\s(\bar{L},\overline{\ad}_k)$ by
\Beq
\bar{\delta}(\bar{x},\bar{y})=\overline{\delta(x,y)}, ~\all \bar{x},\bar{y}\in \bar{L}.
\Eeq
If $\delta_1$ and $\delta_2$ satisfy $\bar{\delta}_1=\bar{\delta}_2$, then $\hat{\delta}\set \delta_1-\delta_2$ is a skew-symmetric biderivation with $\hat{\delta}(L,L)\subset Z(L)$, so $\hat{\delta}$ is a central biderivation on $(L,\alpha)$. Then we derive the following lemma.
\blem
Let $(L,\alpha)$ be a Hom-Lie algebra with $\alpha$ surjective. Up to $\CBider_\s(L,\ad_k)$, the map $\delta\rightarrow \bar{\delta}$ is a $1$-$1$ map from $\Bider_\s(L,\ad_k)$ to $\Bider_\s(\bar{L},\overline{\ad}_k)$.
\elem

Hence, $\Bider_\s(L,\ad_k)$ are determined by $\CBider_\s(L,\ad_k)$ and $\Bider_\s(\bar{L},\overline{\ad}_k)$. We define a sequence of quotient Hom-Lie algebras
\beq\label{seq}
(L^{(0)}, \alpha^{(0)})=(L,\alpha),\cdots, (L^{(r+1)},\alpha^{(r+1)})=(L^{(r)}/Z(L^{(r)}),\overline{\alpha^{(r)}}),
\eeq
and write $(L^{(r)},\ad_k^{(r)},\alpha^{(r)})$ for the quotient adjoint module of $(L^{(r)}, \alpha^{(r)})$.
If there exists $r\in\N$ such that $Z(L^{(r)})=\{0\}$, then $\Bider_\s(L,\ad_k)$ are characterized by $\CBider_\s(L^{(i)},\ad_k^{(i)})(0\leq i\leq r-1)$ and $\Bider_\s(L^{(r)}, \ad_k^{(r)})$, where $\Bider_\s(L^{(r)}, \ad_k^{(r)})$ can be determined by Theorem \ref{BiderofsimpleHomLiealg} in case $(L^{(r)},\alpha^{(r)})$ is perfect and $\alpha$ is invertible.

However, we still don't know how to determine $\Bider_\s(L^{(r)}, \ad_k^{(r)})$ if $(L^{(r)},\alpha^{(r)})$ is not perfect, which inspires us to consider skew-symmetric biderivations on derived Hom-Lie algebras.

Now let $(L,\alpha)$ be an arbitrary centerless Hom-Lie algebra such that $\alpha$ is invertible and $\delta\in\Bider_\s(L,\ad_k)$. Then $\delta$ satisfies
\beq\label{d(a(u),[xy])=[a(x)d(u,y)-[a(y)d(u,x)]}
\delta(\alpha(u),[x,y])=[\alpha^{k+1}(x),\delta(u,y)] -[\alpha^{k+1}(y),\delta(u,x)]\in L', ~\all x,y,u\in L.
\eeq
Thus we could restrict $\delta$ to $L'\times L'$ to get a skew-symmetric biderivation $\delta'\in\Bider_\s(L',\ad_k)$, which makes sense since $L'$ is $\alpha$-invariant.
If $\delta_1,\ \delta_2\in\Bider_\s(L,\ad_k)$ such that $\delta_1|_{L'\times L'}= \delta_2|_{L'\times L'}$, then $\tilde{\delta}\set \delta_1-\delta_2\in\SBider_\s(L,\ad_k)$.

In fact, it suffices to prove $\tilde{\delta}(L,L)\subset Z_L(L')$, since $\tilde{\delta}(L',L')=0$ is obvious. Note that $\tilde{\delta}(L,L')\subset Z_L(L')$ holds by taking $u,y\in L'$ in (\ref{d(a(u),[xy])=[a(x)d(u,y)-[a(y)d(u,x)]}). Then for any $x,y,z,u,v\in L$, we have
\begin{align*}
[[[x,y],\tilde{\delta}(u,v)],\alpha(z)]
&=[[[x,y],z],\alpha\tilde{\delta}(u,v)]+[\alpha([x,y]),[\tilde{\delta}(u,v),z]]\\
&=[[[x,y],z],\tilde{\delta}(\alpha(u),\alpha(v))]-[[\alpha(x),\alpha(y)],[z,\tilde{\delta}(u,v)]]\\
\text{($\alpha$ is isomorphic)}&=[\alpha^k([[\tilde{x},\tilde{y}],\tilde{z}]), \tilde{\delta}(\alpha(u),\alpha(v))] \\ &\relphantom{=}-[[\alpha^{k+1}(\tilde{x}),\alpha^{k+1}(\tilde{y})], [\alpha^k(\tilde{z}),\tilde{\delta}(u,v)]]\\
\text{(by Lemmas \ref{xyzw=0} and \ref{d(u,[xy])-ud(x,y)})}&=-[[\alpha^{k+1}(u),\alpha^{k+1}(v)], \tilde{\delta}([\tilde{x},\tilde{y}],\tilde{z})] \\&\relphantom{=}-[[\alpha^{k+1}(\tilde{x}),\alpha^{k+1}(\tilde{y})],\tilde{\delta}(\tilde{z},[u,v])]\\
(\text{by $\tilde{\delta}(L,L')\subset Z_L(L')$})&=0.
\end{align*}
It follows that $\tilde{\delta}(L,L)\subset Z_L(L')$, since $(L,\alpha)$ is centerless. Thus, $\tilde{\delta}\in\SBider_\s(L,\ad_k)$, and we deduce the following lemma.

\blem
Suppose that $(L,\alpha)$ is a centerless Hom-Lie algebra with $\alpha$ invertible. Then any $\delta\in\Bider_\s(L,\ad_k)$, up to $\SBider_\s(L,\ad_k)$, is an extension of a unique $\delta'\in\Bider_\s(L',\ad_k)$.
\elem

Hence, for centerless but not perfect Hom-Lie algebra $(L^{(r)}, \alpha^{(r)})$, $\Bider_\s(L^{(r)}, \ad_k^{(r)})$ are determined by $\SBider_\s(L^{(r)}, \ad_k^{(r)})$ and $\Bider_\s({L^{(r)}}',\ad_k^{(r)})$. Taking $(L,\alpha)=({L^{(r)}}', \alpha^{(r)})$ and repeating the above arguments based on (\ref{seq}), then we continue the algorithms.

Now we apply our method to several concrete examples.

\beg
Let $(L,\alpha)$ belong to one class of multiplicative Heisenberg Hom-Lie algebras in \cite[Corollary 2.3]{Alvarez&Cartes}. Specifically, $L=\F e_1\oplus \F e_2\oplus \F e_3$ as a vector space with $[e_1,e_2]=e_3$, $[e_2,e_3]=[e_1,e_3]=0$, and
$\alpha=\left(
 \begin{array}{ccc}
 \lambda & 0 & 0 \\
 1 & \lambda & 0 \\
 0 & 0 & \lambda^2 \\
 \end{array}
 \right)$
for some nonzero $\lambda\in\F$.
Note that $Z(L)=\F e_3$. As in (\ref{seq}), we have
\Beq
(L^{(0)}, \alpha^{(0)})=(L,\alpha),\quad
(L^{(1)}, \alpha^{(1)})=(L/Z(L),\bar{\alpha})= \left(\F \bar{e}_1\oplus \F \bar{e}_2, \left( \begin{array}{cc}
\lambda & 0 \\
1 & \lambda \\
\end{array}
\right)
\right).
\Eeq
Since $(L^{(1)}, \alpha^{(1)})$ is abelian,
$(L^{(2)},\alpha^{(2)})=(\{0\},0).$ Then $\Bider_\s(L^{(2)},\ad_k^{(2)})=\{0\}$, and we need $\CBider_\s(L^{(0)},\ad_k^{(0)})$ and $\CBider_\s(L^{(1)},\ad_k^{(1)})$ to determine $\Bider_\s(L,\ad_k)$.

Let $\delta^{(1)}\in\CBider_\s(L^{(1)},\ad_k^{(1)})$. Assume $\delta^{(1)}(\bar{e}_1, \bar{e}_2)=k_1\bar{e}_1+k_2\bar{e}_2.$ Note that $\delta^{(1)}$ satisfies
\Beq
\alpha^{(1)}\delta^{(1)}(\bar{e}_1, \bar{e}_2)=\delta^{(1)}(\alpha^{(1)}(\bar{e}_1), \alpha^{(1)}(\bar{e}_2)),
\Eeq
which implies
\Beq
k_1\lambda\bar{e}_1+(k_1+k_2\lambda)\bar{e}_2 =\lambda^2(k_1\bar{e}_1+k_2\bar{e}_2).
\Eeq
Then
\beq\label{tbiderofL1}
\delta^{(1)}(\bar{e}_1, \bar{e}_2)
=\left\{
\begin{array}{ll}
k_2\bar{e}_2, & \hbox{$\lambda=1$;} \\
0, & \hbox{$\lambda\neq1$.}
\end{array}
\right.
\eeq
The above $\delta^{(1)}$ clearly satisfies (\ref{fixleftder}) or (\ref{fixrightder}) since $(L^{(1)}, \alpha^{(1)})$ is an abelian Hom-Lie algebra. Note also that $\Bider_\s(L^{(2)},\ad_k^{(2)})=\{0\}$. Then any $\delta\in\Bider_\s(L^{(1)},\ad_k^{(1)})$ must be some $\delta^{(1)}$ in (\ref{tbiderofL1}).

Let $\delta^{(0)}\in\CBider_\s(L^{(0)},\ad_k^{(0)})$.
By Definition \ref{CBider&SBider}, Remark \ref{CBider} and $L'=Z(L)=\F e_3$, we could set
\beq\label{tbiderofL}
\delta^{(0)}(e_1,e_2)=k_3e_3, ~~\delta^{(0)}(e_1,e_3)=\delta^{(0)}(e_2,e_3)=0.
\eeq
It is straightforward to verify that the above $\delta^{(0)}$ satisfies (\ref{bd=da})-(\ref{fixrightder}).

Now let $\hat{\delta}\in\Bider_\s(L,\ad_k)$. Then $\hat{\delta}$ induces a skew-symmetric biderivation that belongs to $\Bider_\s(L^{(1)},\ad_k^{(1)})$, which is of the form in (\ref{tbiderofL1}). Set
\beq\label{biderofL}
\hat{\delta}(e_1, e_2)
=\left\{
\begin{array}{ll}
k_2e_2, & \hbox{$\lambda=1$;} \\
0, & \hbox{$\lambda\neq1$,}
\end{array}
\right. ~~ \hat{\delta}(e_1,e_3)=\hat{\delta}(e_2,e_3)=0 ~(\forall~\lambda\in\F^*).
\eeq
One could check that $\hat{\delta}$ satisfies (\ref{bd=da})-(\ref{fixrightder}).

Therefore, every $\delta\in\Bider_\s(L,\ad_k)$ is of the form \Beq
\delta=\delta^{(0)}+\hat{\delta},
\Eeq
for some $\delta^{(0)}$ in (\ref{tbiderofL}) and $\hat{\delta}$ in (\ref{biderofL}).
In particular,
\Beq
\delta(e_1, e_2)
=\left\{
\begin{array}{ll}
k_2e_2+k_3e_3, & \hbox{$\lambda=1$;} \\
k_3e_3, & \hbox{$\lambda\neq1$,}
\end{array}
\right. ~~ \delta(e_1,e_3)=\delta(e_2,e_3)=0 ~(\forall~\lambda\in\F^*).
\Eeq
\eeg

\beg\label{exofHomLiealg}
Suppose that $(L,\alpha)$ belongs to the second class of multiplicative Hom-Lie algebras in \cite[Theorem 4.7]{Chen&Zhang}.
For simplicity of calculations, we consider $L=\F x\oplus \F y\oplus \F z$ as a vector space with $[x,y]=y$, $[x,z]=[y,z]=0$ and
$\alpha=\left(
 \begin{array}{ccc}
 1 & 0 & 0 \\
 a & \lambda & 0 \\
 b & 0 & \mu \\
 \end{array}
 \right)$
for $a,b\in\F$ and $\lambda,\mu\in\F^*$. Note that $Z(L)=\F z$. As in (\ref{seq}), we have
\Beq
(L^{(0)}, \alpha^{(0)})=(L,\alpha),\quad
(L^{(1)}, \alpha^{(1)})=\left(\F \bar{x}\oplus \F \bar{y},
\left(
\begin{array}{cc}
1& 0 \\
a & \lambda \\
\end{array}
\right)
\right).
\Eeq
Then $(L^{(1)}, \alpha^{(1)})$ is centerless but not perfect since ${L^{(1)}}'=\F \bar{y}$. Hence $\Bider_\s(L^{(1)},\ad_k^{(1)})$ is determined by $\SBider_\s(L^{(1)},\ad_k^{(1)})$ and $\Bider_\s({L^{(1)}}',\ad_k^{(1)})$, where the latter is $\{0\}$ since dim${L^{(1)}}'=1$. So we only need $\CBider_\s(L^{(0)},\ad_k^{(0)})$ and $\SBider_\s(L^{(1)},\ad_k^{(1)})$ to determine $\Bider_\s(L,\ad_k)$.

Let $\tilde{\delta}\in\SBider_\s(L^{(1)},\ad_k^{(1)})$. Note that $Z_{L^{(1)}}({L^{(1)}}')={L^{(1)}}'=\F \bar{y}$. We set
\beq\label{sbider}
\tilde{\delta}(\bar{x},\bar{y})=k\bar{y},
\eeq
which can be verified to be a skew-symmetric biderivation. Then all skew-symmetric biderivations in $\Bider_\s(L^{(1)},\ad_k^{(1)})$ are of the form given by (\ref{sbider}), since there is no nonzero element in $\Bider_\s({L^{(1)}}',\ad_k^{(1)})$.

Let $\delta^{(0)}\in\CBider_\s(L^{(0)},\ad_k^{(0)})$. Since $L'=\F y$ and $Z(L)=\F z$,
we could set
\Beq
\delta^{(0)}(x,z)=lz, ~~\delta^{(0)}(x,y)=\delta^{(0)}(y,z)=0.
\Eeq
It is easy to check that the above $\delta^{(0)}$ is a skew-symmetric biderivation, then such $\delta^{(0)}$ are precisely all elements in $\CBider_\s(L^{(0)},\ad_k^{(0)})$.

Set
\beq\label{biderof2ndL}
\delta(x, y)
=ky, ~~ \delta(x,z)=\delta(y,z)=0.
\eeq
Then $\delta\in\Bider_\s(L,\ad_k)$ and $\delta$ induces a skew-symmetric biderivation in $\Bider_\s(L^{(1)},\ad_k^{(1)})$, which is of the form in (\ref{sbider}).

Therefore, every $\delta\in\Bider_\s(L,\ad_k)$ is of the following form
\Beq
\delta(x,y)=ky, ~~\delta(x,z)=lz, ~~\delta(y,z)=0.
\Eeq
\eeg

\beg
Consider the loop algebra $L=\sl_2\otimes \F[t, t^{-1}]$ with $[x\otimes t^m, y\otimes t^n]=[x, y]\otimes t^{m+n}$ and $\alpha=\check{\alpha}\otimes \id$, where $\check{\alpha}$ is an involution of $\sl_2$ such that
\Beq
\check{\alpha}(e)=-e, ~~\check{\alpha}(f)=-f, ~~\check{\alpha}(h)=h.
\Eeq
Then $(L,\alpha)$ is an infinite-dimensional Hom-Lie algebra (cf. \cite[Corollaries 3.1 and 3.3]{Xie&Liu}). Moreover, $(L,\alpha)$ is perfect and centerless, so each $\delta\in\Bider_\s(L,\ad_k)$ is of the form
\Beq
\delta(x,y)=\alpha^{-1}\gamma([x,y]), ~\all x, y\in L,
\Eeq
for some $\gamma\in \Cent (L, \ad_k)$.

Let $\gamma\in \Cent (L, \ad_k)$. Then $\gamma\alpha=\alpha\gamma$ and for all $a, b\in\sl_2$ and $m, n\in\Z$,
\Beq
\gamma([a\otimes t^m, b\otimes t^n])=[\alpha^{k+1}(a\otimes t^m),\gamma(b\otimes t^n)]
\Eeq
Take $a\otimes t^m=h\otimes 1$. It follows that $\gamma([h,b]\otimes t^n)=[h\otimes 1, \gamma(b\otimes t^n)]$. In particular,
\Beq
-2\gamma(f\otimes t^n)=[h\otimes 1, \gamma(f\otimes t^n)], ~~
0=[h\otimes 1, \gamma(h\otimes t^n)], ~~
2\gamma(e\otimes t^n)=[h\otimes 1, \gamma(e\otimes t^n)],
\Eeq
which implies that $\gamma(f\otimes t^n)$, $\gamma(h\otimes t^n)$ and $\gamma(e\otimes t^n)$ are eigenvectors for the operator $[h\otimes 1, -]$ with eigenvalues $-2$, $0$ and $2$
respectively. Then there exist $\Phi_f^n(t), \Phi_h^n(t), \Phi_e^n(t)\in \F[t, t^{-1}]$ such that
\Beq
\gamma(f\otimes t^n)=f\otimes \Phi_f^n(t), ~~
\gamma(h\otimes t^n)=h\otimes \Phi_h^n(t), ~~
\gamma(e\otimes t^n)=e\otimes \Phi_e^n(t).
\Eeq
Note that
\begin{align*}
-2f\otimes \Phi_f^n(t)
&=-2\gamma(f\otimes t^n)=\gamma([h\otimes t^n, f\otimes 1])\\
&=[h\otimes t^n,\gamma(f\otimes 1)]=[h\otimes t^n,f\otimes \Phi_f^0(t)]=-2f\otimes t^n\Phi_f^0(t).
\end{align*}
Then $\Phi_f^n(t)=t^n\Phi_f^0(t), \all n\in\Z$. Similarly, $\Phi_e^n(t)=t^n\Phi_e^0(t)$. Note also that
\begin{align*}
2f\otimes \Phi_f^n(t)
&=2\gamma(f\otimes t^n)=\gamma([f\otimes 1, h\otimes t^n])\\
&=[(-1)^{k+1}f\otimes 1, \gamma(h\otimes t^n)]=[(-1)^{k+1}f\otimes 1, h\otimes \Phi_h^n(t)]=2f\otimes (-1)^{k+1}\Phi_h^n(t).
\end{align*}
Then we have $\Phi_h^n(t)=t^n\Phi_h^0(t)$ and $\Phi_f^0(t)=(-1)^{k+1}\Phi_h^0(t)$. Substituting $f$ with $e$ in the above equality one obtains $\Phi_e^0(t)=(-1)^{k+1}\Phi_h^0(t)$. Hence
\Beq
(-1)^{k+1}\Phi_f^0(t)=(-1)^{k+1}\Phi_e^0(t)=\Phi_h^0(t).
\Eeq
Set $\Phi(t)=\Phi_h^0(t)$ and define $\phi\in\End(\F[t,t^{-1}])$ by $\phi(g(t))=\Phi(t)g(t)$.
It follows that $\gamma=\check{\alpha}^{k+1}\otimes \phi$.

Conversely, for any $\Phi(t)\in\F[t,t^{-1}]$, the $\gamma$ defined as above belongs to $\Cent (L, \ad_k)$. Therefore,
\Beq
\Cent (L, \ad_k)=\{\check{\alpha}^{k+1}\otimes \phi ~|~ \Phi(t)\in\F[t,t^{-1}]\},
\Eeq
and so every $\delta\in\Bider_\s(L,\ad_k)$ is of the following form:
\Beq
\delta(x,y)=\alpha^{-1}(\check{\alpha}^{k+1}\otimes \phi)([x,y])=(\check{\alpha}^{-1}\otimes \id)(\check{\alpha}^{k+1}\otimes \phi)([x,y])=(\check{\alpha}^{k}\otimes \phi) ([x,y]), ~\all x, y\in L.
\Eeq
\eeg
\section{Commuting linear maps on Hom-Lie algebras}
In this section, we consider linear maps from a Hom-Lie algebra $(L,\alpha)$ to an $(L,\alpha)$-module $(V, \rho, \beta)$ that are also closely related to $\Cent (L, V)$.
\bdefn\label{comm map}
Let $(L,\alpha)$ be a Hom-Lie algebra and $(V, \rho, \beta)$ an $(L,\alpha)$-module. A linear map $f: L\rightarrow V$ is called a \textit{commuting linear map} if
\Beq
\alpha(x) f(x)=0 ~~{\rm and}~~ \beta\circ f=f\circ \alpha, ~\all x\in L.
\Eeq
\edefn
Denote by $\Com(L, V)$ the set of all commuting linear maps $f: L\rightarrow V$. In particular, we write $\Com(L, \ad_k)$ as the set of all commuting linear maps on $(L, \alpha)$ for the adjoint module $(L, \ad_k, \alpha)$,

Clearly, $\Cent (L, V)\subseteq \Com(L, V)$. We will show $\Cent (L, V)=\Com(L, V)$, under the assumption that $\beta$ is invertible and $Z_V(L')={0}$.
\blem\label{[a(v)a(w)]a(u)(f([xy])-a(x)f(y))=0}
Let $(L,\alpha)$ be a Hom-Lie algebra and $(V, \rho, \beta)$ an $(L,\alpha)$-module with $\beta$ invertible. If $f\in\Com(L, V)$, then
\Beq
[\alpha(v),\alpha(w)]\alpha^2(u)(f([x,y])-\alpha(x)f(y))=0, ~\all x,y,z,u,v\in L.
\Eeq
\elem

\bpf
By the definition of the commuting linear map, we get
\beq\label{xf(y)=-yf(x)}
\alpha(x) f(y)=-\alpha(y) f(x), ~\all x,y\in L,
\eeq
which induces a skew-symmetric bilinear map $\delta: L\times L\rightarrow V$ defined by
\beq \label{delta(x,y)=betaalpha(x)f(y)}
\delta(x,y)=\beta^{-1}(\alpha(x) f(y)), ~\all x,y\in L.
\eeq
Moreover, $\delta$ is a skew-symmetric biderivation, since \begin{align*}
\beta\delta(x,y)&=\alpha(x) f(y)
=\alpha(x)\beta^{-1} \beta f(y)\\
\text{(by Definition \ref{comm map})}&=\alpha(x)\beta^{-1} f(\alpha(y))\\
\text{(by (\ref{bp=pab}))}&=\beta^{-1}\alpha^2(x) f(\alpha(y))\\ &=\delta(\alpha(x),\alpha(y))
\end{align*}
and
\begin{align*}
\delta(\alpha(x),[y,z])&=\beta^{-1}\big(\alpha^2(x)f([y,z])\big)\\
\text{(by (\ref{xf(y)=-yf(x)}))}&=-\beta^{-1}\big(\alpha([y,z])f(\alpha(x))\big) =-\beta^{-1}\big([\alpha(y),\alpha(z)]\beta f(x)\big)\\
\text{(by (\ref{pb=pap-pap}))}&=-\beta^{-1}\big(\alpha^2(y)\alpha(z)f(x)-\alpha^2(z)\alpha(y)f(x)\big)\\
\text{(by (\ref{bp=pab}))}&=-\alpha(y)\beta^{-1}\big(\alpha(z)f(x)\big) +\alpha(z)\beta^{-1}\big(\alpha(y)f(x)\big)\\
&=-\alpha(y)\delta(z,x)+\alpha(z)\delta(y,x)\\
&=\alpha(y)\delta(x,z)-\alpha(z)\delta(x,y).
\end{align*}
Note that Lemma \ref{d(u,[xy])-ud(x,y)} says that for any $v, w\in L$, $[v,w](\delta(\alpha(u),[x,y])-\alpha(u)\delta(x,y))=0$. It follows that for all $x,y,u,v,w\in L$,
\begin{align*}
0&=\beta\big([v,w](\delta(\alpha(u),[x,y])-\alpha(u)\delta(x,y))\big)\\
\text{(by (\ref{bp=pab}))}&=\alpha([v,w])\beta\big(\delta(\alpha(u),[x,y])-\alpha(u)\delta(x,y)\big)\\
\text{(by (\ref{bp=pab}))}&=\alpha([v,w])\big(\beta\delta(\alpha(u),[x,y])-\alpha^2(u)\beta\delta(x,y)\big)\\
(\text{by (\ref{delta(x,y)=betaalpha(x)f(y)})})&=[\alpha(v),\alpha(w)]\big(\alpha^2(u)f([x,y])-\alpha^2(u)\alpha(x)f(y)\big)\\
&=[\alpha(v),\alpha(w)]\alpha^2(u)\big(f([x,y])-\alpha(x)f(y)\big),
\end{align*}
as required.
\epf

\bthm\label{ComofHomLiealg}
Let $(L,\alpha)$ be a Hom-Lie algebra with $\alpha$ surjective and $(V, \rho, \beta)$ an $(L,\alpha)$-module with $\beta$ invertible. If $Z_V(L')={0}$, then $\Cent (L, V)=\Com(L, V)$.
\ethm
\bpf
It suffices to prove $\Com(L, V)\subseteq\Cent (L, V)$.

Take $f\in\Com(L, V)$.
Note that $Z_V(L)\subseteq Z_V(L')={0}$. Then $f([x,y])=\alpha(x)f(y)$ by Lemma \ref{[a(v)a(w)]a(u)(f([xy])-a(x)f(y))=0}, so $f\in \Cent (L, V)$.
\epf

Now we would like to describe $\Com(L, V)$, which need the following notions of central and special commuting linear maps.

\bdefn
Let $(L,\alpha)$ be a Hom-Lie algebra and $(V, \rho, \beta)$ an $(L,\alpha)$-module. A commuting linear map $f\in\Com(L, V)$ is called \textit{central} if
\Beq
f(L)\subset Z_V(L);
\Eeq
$f$ is called \textit{special} if
\Beq
f(L)\subset Z_V(L')~~\text{and}~~ f(L')=0.
\Eeq
Denote by $\CCom(L, V)$ and $\SCom(L, V)$ the sets consisting of all central and special commuting linear maps in $\Com(L, V)$, respectively.

Clearly, every linear map $f: L\rightarrow Z_V(L)$ satisfying $\beta\circ f=f\circ \alpha$ is automatically a central commuting linear map.
\edefn
Note that $Z_V(S)$ is a submodule of $V$ when $S$ is an ideal of $(L, \alpha)$ such that $\alpha(S)=S$. In fact, for any $s\in S$ and $v\in Z_V(S)$, there exists $t\in S$ satisfying $s=\alpha(t)$. Then $Z_V(S)$ is $\beta$-invariant by
\Beq
s\beta(v)=\alpha(t)\beta(v)=\beta(tv)=0;
\Eeq
$Z_V(S)$ is $L$-invariant by
\Beq
sxv=\alpha(t)xv=[t,x]\beta(v)+\alpha(x)tv=0,~\all x\in L.
\Eeq

Now consider a Hom-Lie algebra $(L,\alpha)$ with $\alpha$ surjective. Since $L'$ is an ideal of $(L,\alpha)$, $Z_V(L')$ becomes a submodule of $(V, \rho, \beta)$, which induces a quotient module $(V/Z_V(L'),\bar{\rho}, \bar{\beta})$. Then for any $f\in\Com(L, V)$, we can define
$\bar{f}\in\Com(L, V/Z_V(L'))$ by $\bar{f}(x)= f(x)+Z_V(L')$.

If $f_1,\ f_2\in\Com(L, V)$ such that $\bar{f_1}=\bar{f_2}$, then $f\set f_1-f_2$ satisfies $f(L)\subset Z_V(L')$. Note that
\Beq
0=\alpha([x,y])f(z)=-\alpha(z)f([x,y]), ~\all x,y,z\in L,
\Eeq
which implies $f(L')\subset Z_V(L)$. Let $\hat{f}: L\rightarrow Z_V(L)\in\CCom(L, V)$ such that $\hat{f}|_{L'}=f|_{L'}$. Then $\tilde{f}\set f-\hat{f}$ satisfies
\Beq
\tilde{f}(L')=0 ~~\text{and}~~ \tilde{f}(L)\subset Z_V(L'),
\Eeq
that is, $\tilde{f}\in\SCom(L, V)$. Hence $f=f_1-f_2=\hat{f}+\tilde{f}\in\CCom(L, V)+\SCom(L, V)$. Then we prove the following property.
\bprop
Let $(L,\alpha)$ be a Hom-Lie algebra with $\alpha$ surjective and $(V, \rho, \beta)$ an $(L,\alpha)$-module. Up to $\CCom(L, V)+\SCom(L, V)$, the map $f\rightarrow \bar{f}$ is a $1$-$1$ map from $\Com(L, V)$ to $\Com(L, V/Z_V(L'))$.
\eprop

Hence, in case $(V, \rho, \beta)$ is an $(L,\alpha)$-module with $\alpha$ surjective, $\Com(L, V)$ are determined by $\CCom(L, V)$, $\SCom(L, V)$ and $\Com(L, V/Z_V(L'))$. Then we could give an algorithm for describing $\Com(L, V)$ as follows.
Define a sequence of quotient modules
\Beq
(V^{[0]},\rho^{[0]},\beta^{[0]})=(V,\rho,\beta), \cdots, (V^{[r+1]},\rho^{[r+1]},\beta^{[r+1]}) =(V^{[r]}/Z_{V^{[r]}}(L'),\overline{\rho^{[r]}},\overline{\beta^{[r]}}).
\Eeq
If there exists $r\in\N$ such that $Z_{V^{[r]}}(L')=\{0\}$, then $\Com(L, V)$ are determined by $\CCom(L, V^{(i)})$, $\SCom(L, V^{(i)})$ ($0\leq i\leq r-1$) and $\Com(L, V^{[r]})$, where $\Com(L, V^{[r]})$ is precisely $\Cent (L, V^{[r]}$) in case $\beta$ is required to be invertible.

Now we apply our method to a concrete example.
\beg
Let $(L,\alpha)$ be the multiplicative Hom-Lie algebras in Example \ref{exofHomLiealg} with
$\alpha=\left(
 \begin{array}{ccc}
 1 & 0 & 0 \\
 0 & \lambda & 0 \\
 0 & 0 & \mu \\
 \end{array}
 \right)$
and $(V, \rho, \beta)=(L, \ad_k, \alpha)$. Let's determine $\Com(L, \ad_k)$.

Note that $L'=\F y$ and $Z_L(L')=\F y\oplus\F z$. Then $(L^{[1]}, \ad_k^{[1]}, \alpha^{[1]})=(\F \bar{x}, 0, \id)$. For $f^{[1]}\in\Com(L, L^{[1]})$, suppose
\beq\label{Com(L1)}
f^{[1]}(x)=k_1\bar{x},~~f^{[1]}(y)=k_2\bar{x},~~f^{[1]}(z)=k_3\bar{x}.
\eeq
Then $f^{[1]}$ clearly satisfies $[\overline{\alpha^{k+1}(x)}, f^{[1]}(x)]=[\overline{\alpha^{k+1}(y)}, f^{[1]}(y)]=[\overline{\alpha^{k+1}(z)}, f^{[1]}(z)]=\bar{0}$; $f^{[1]}$ satisfies $\id\circ f^{[1]}=f^{[1]}\circ\alpha$, iff
\beq\label{coe}
k_2(\lambda-1)=k_3(\mu-1)=0.
\eeq
Hence every $f^{[1]}\in\Com(L, L^{[1]})$ is defined by (\ref{Com(L1)}) and (\ref{coe}).

Suppose that $f\in\Com(L, \ad_k)$ such that $\bar{f}=f^{[1]}$. Set $f(y)=k_2x+k_2'y+k_2''z$. Then $k_2=0$, which follows from \Beq
0=[\alpha^{k+1}(y), f(y)]=[\lambda^{k+1}y,k_2x+k_2'y+k_2''z]=-k_2\lambda^{k+1}y.
\Eeq
Hence $f^{[1]}\in\Com(L, L^{[1]})$ is induced from some $f\in\Com(L, \ad_k)$ only if $f^{[1]}\in\Com(L, L^{[1]})$ is defined by (\ref{Com(L1)}) such that $k_3(\mu-1)=0$. Define a linear map $f: L\rightarrow L$ by
\Beq
f(x)=k_1x, ~~f(y)=0, ~~f(z)=k_3x, ~~(k_1, k_3\in \F, k_3(\mu-1)=0).
\Eeq
It is straightforward to verify that $f\in\Com(L, \ad_k)$.

Now it remains to compute $\CCom(L, \ad_k)$ and $\SCom(L, \ad_k)$. Let $\hat{f}\in\CCom(L, \ad_k)$ and $\tilde{f}\in\SCom(L, \ad_k)$. Since $Z_L(L)=\F z$, $L'=\F y$ and $Z_L(L')=\F y\oplus\F z$, set
\Beq
\hat{f}(x)=l_1z,~~\hat{f}(y)=l_2z,~~\hat{f}(z)=l_3z;
\Eeq
\Beq \tilde{f}(x)=a_{11}y+a_{21}z,~~\tilde{f}(y)=0,~~\tilde{f}(z)=a_{12}y+a_{22}z.
\Eeq
By a direct computation, $\hat{f}$ is a commuting linear map iff $l_1(\mu-1)=l_2(\lambda-\mu)=0$; $\tilde{f}$ satisfies $[\alpha^{k+1}(x),\tilde{f}(x)] =[\alpha^{k+1}(y),\tilde{f}(y)] =[\alpha^{k+1}(z),\tilde{f}(z)]=0$ iff $a_{11}=0$; $\tilde{f}$ satisfies $\alpha\circ \tilde{f}=\tilde{f}\circ\alpha$, iff
\Beq
a_{11}(\lambda-1)=a_{21}(\mu-1)=a_{12}(\lambda-\mu)=0.
\Eeq
Hence, $\tilde{f}$ is of the form
\Beq
\tilde{f}(x)=a_{21}z, ~~\tilde{f}(y)=0, ~~\tilde{f}(z)=a_{12}y+a_{22}z,
\Eeq
where $a_{21}, a_{12}, a_{22}\in\F$ such that $a_{21}(\mu-1)=a_{12}(\lambda-\mu)=0$.

Therefore, every $f\in\Com(L, \ad_k)$ is of the form
\Beq
f(x)=c_1x+c_2z,~~f(y)=c_3z,~~f(z)=c_4x+c_5y+c_6z,
\Eeq
for all $c_1, c_2, c_3, c_4, c_5, c_6\in \F$ such that
\Beq
c_2(\mu-1)=c_3(\lambda-\mu)=c_4(\mu-1)=c_5(\lambda-\mu)=0.
\Eeq
\eeg

Finally, we would like to describe the relationship between skew-symmetric biderivations and commuting linear maps on a Hom-Lie algebra $(L, \alpha)$ for the adjoint module $(L, \ad_k, \alpha)$.

\bprop
Let $(L,\alpha)$ be a Hom-Lie algebra with $\alpha$ invertible and $(L, \ad_k, \alpha)$ an adjoint $(L,\alpha)$-module such that every $\delta\in\Bider_\s(L,\ad_k)$ is of the form $\delta(x,y)=\alpha^{-1}\gamma([x,y])$, where $\gamma\in \Cent (L, \ad_k)$. Then every $f\in\Com(L,\ad_k)$ is of the form $f=\gamma+\mu$ for some $\gamma\in \Cent (L, \ad_k)$ and $\mu\in\CCom(L, \ad_k)$.
\eprop

\bpf
Define
\Beq \delta(x,y)=\alpha^{-1}[\alpha^{k+1}(x),f(y)],
\Eeq
which belongs to $\Bider_\s(L,\ad_k)$, by the proof of Lemma \ref{[a(v)a(w)]a(u)(f([xy])-a(x)f(y))=0} when $(V,\rho,\beta)$ is taken as $(L, \ad_k, \alpha)$.

Thus there exists $\gamma\in \Cent (L, \ad_k)$ such that
$\gamma([x,y])=\alpha\delta(x,y)=[\alpha^{k+1}(x),f(y)]$. Note that $\gamma([x,y])=[\alpha^{k+1}(x),\gamma(y)]$, which implies $[\alpha^{k+1}(x),(f-\gamma)(y)]=0$, then $f-\gamma\in\CCom(L, \ad_k)$.
\epf


\begin{thebibliography}{99}
\bibitem{Agrebaoui&Benali&Makhlouf} B. Agrebaoui, K. Benali, A. Makhlouf, Representations of simple Hom-Lie algebras. \textit{J. Lie Theory} 29 (2019), no. 4, 1119--1135.
\bibitem{Aizawa&Sato} N. Aizawa, H. Sato, $q$-deformation of the Virasoro algebra with central extension. \textit{Phys. Lett. B} 256 (1991), no. 2, 185--190.
\bibitem{Alvarez&Cartes} M. Alvarez, F. Cartes, Cohomology and deformations for the Heisenberg Hom-Lie algebras. \textit{Linear Multilinear Algebra} 67 (2019), no. 11, 2209--2229.
\bibitem{Bresar} M. Bre\v{s}ar, Commuting maps: a survey. \textit{Taiwanese J. Math.} 8 (2004), no. 3, 361--397.
\bibitem{Bresar&Martindale&Miers} M. Bre\v{s}ar, W. Martindale III, C. Miers, Centralizing maps in prime rings with involution. \textit{J. Algebra} 161 (1993), no. 2, 342--357.
\bibitem{Bresar&Zhao} M. Bre\v{s}ar, K. Zhao, Biderivations and commuting linear maps on Lie algebras. \textit{J. Lie Theory} 28 (2018), no. 3, 885--900.
\bibitem{Casas&GarciaMartinez} J. Casas, X. Garc\'{\i}a-Mart\'{\i}nez, Abelian extensions and crossed modules of Hom-Lie algebras. \textit{J. Pure Appl. Algebra} 224 (2020), no. 3, 987--1008.
\bibitem{Chaichian&Kulish} M. Chaichian, P. Kulish, J. Lukierski, $q$-deformed Jacobi identity, $q$-oscillators and $q$-deformed infinite-dimensional algebras. \textit{Phys. Lett. B} 237 (1990), no. 3--4, 401--406.
\bibitem{Chen&Han} X. Chen, W. Han, Classification of multiplicative simple Hom-Lie algebras. \textit{J. Lie Theory} 26 (2016), no. 3, 767--775.
\bibitem{Chen&Zhang} Y. Chen, R. Zhang, A commutative algebra approach to multiplicative Hom-Lie algebras. arXiv:1907.02415v3.
\bibitem{Chen} Z. Chen, Biderivations and linear commuting maps on simple generalized Witt algebras over a field. \textit{Electron J. Linear Algebra} 31 (2016), 1--12.
\bibitem{Curtright&Zachos} T. Curtright, C. Zachos, Deforming maps for quantum algebras. \textit{Phys. Lett. B} 243 (1990), no. 3, 237--244.
\bibitem{Farkas&Letzter} D. Farkas, G. Letzter, Ring theory from symplectic geometry. \textit{J. Pure Appl. Algebra} 125 (1998), no. 1--3, 155--190.
\bibitem{Han&Wang&Xia} X. Han, D. Wang, C. Xia, Linear commuting maps and biderivations on the Lie algebras $W(a,b)$. \textit{J. Lie Theory} 26 (2016), no. 3, 777--786.
\bibitem{Hartwig&Larsson&Silvestrov} J. Hartwig, D. Larsson, S. Silvestrov, Deformations of Lie algebras using $\sigma$-derivations. \textit{J. Algebra} 295 (2006), no. 2, 314--361.
\bibitem{Hu} N. Hu, $q$-Witt algebras, $q$-Lie algebras, $q$-holomorph structure and representations. \textit{Algebra Colloq.} 6 (1999), no. 1, 51--70.
\bibitem{Jin&Li} Q. Jin, X. Li, Hom-Lie algebra structures on semi-simple Lie algebras. \textit{J. Algebra} 319 (2008), no. 4, 1398--1408.
\bibitem{Kassel} C. Kassel, Cyclic homology of differential operators, the Virasoro algebra and a $q$-analogue. \textit{Comm. Math. Phys.} 146 (1992), no. 2, 343--356.
\bibitem{Larsson&Silvestrov} D. Larsson, S. Silvestrov, Quasi-hom-Lie algebras, central extensions and 2-cocycle-like identities. \textit{J. Algebra} 288 (2005), no. 2, 321--344.
\bibitem{Laurent-Gengoux&Makhlouf&Teles} C. Laurent-Gengoux, A. Makhlouf, J. Teles, Universal algebra of a Hom-Lie algebra and group-like elements. \textit{J. Pure Appl. Algebra} 222 (2018), no. 5, 1139--1163.
\bibitem{Laurent-Gengoux&Teles} C. Laurent-Gengoux, J. Teles, Hom-Lie algebroids. \textit{J. Geom. Phys.} 68 (2013), 69--75.
\bibitem{Liu&Guo&Zhao} X. Liu, X. Guo, K. Zhao, Biderivations of the block Lie algebras. \textit{Linear Algebra Appl.} 538 (2018), 43--55.
\bibitem{Makhlouf&Silvestrov} A. Makhlouf, S. Silvestrov, Notes on 1-parameter formal deformations of Hom-associative and Hom-Lie algebras. \textit{Forum Math.} 22 (2010), no. 4, 715--739.
\bibitem{Makhlouf&Zusmanovich} A. Makhlouf, P. Zusmanovich, Hom-Lie structures on Kac-Moody algebras. \textit{J. Algebra} 515 (2018), 278--297.
\bibitem{Maksa} Gy. Maksa, A remark on symmetric biadditive functions having nonnegative diagonalization. \textit{Glasnik Mat. Ser. III} 15 (1980), no. 2, 279--282.
\bibitem{Peyghan&Nourmohammadifar} E. Peyghan, L. Nourmohammadifar, Para-K\"{a}hler hom-Lie algebras. \textit{J. Algebra Appl.} 18 (2019), no. 3, 1950044, 24 pp.
\bibitem{Sheng} Y. Sheng, Representations of hom-Lie algebras. \textit{Algebr. Represent. Theory} 15 (2012), no. 6, 1081--1098.
\bibitem{Skosyrskii} V. Skosyrski\u{\i}, Strongly prime noncommutative Jordan algebras (Russian). \textit{Trudy Inst. Mat. (Novosibirsk)} 16 (1989), 131--164, 198--199.
\bibitem{Tang} X. Tang, Biderivations of finite-dimensional complex simple Lie algebras. \textit{Linear Multilinear Algebra} 66 (2018), no. 2, 250--259.
\bibitem{Vukman} J. Vukman, Symmetric bi-derivations on prime and semi-prime rings. \textit{Aequationes Math.} 38 (1989), no. 2--3, 245--254.
\bibitem{Wang&Yu} D. Wang, X. Yu, Biderivations and linear commuting maps on the Schr\"{o}dinger-Virasoro Lie algebra. \textit{Comm. Algebra} 41 (2013), no. 6, 2166--2173.
\bibitem{Wang&Yu&Chen} D. Wang, X. Yu, Z. Chen, Biderivations of the parabolic subalgebras of simple Lie algebras. \textit{Comm. Algebra} 39 (2011), no. 11, 4097--4104.
\bibitem{Xie&Liu} W. Xie, W. Liu, Hom-structures on simple graded Lie algebras of finite growth. \textit{J. Algebra Appl.} 16 (2017), no. 8, 1750154, 18 pp.
\bibitem{Yau} D. Yau, Hom-algebras and homology. \textit{J. Lie Theory} 19 (2009), no. 2, 409--421.
\end{thebibliography}
\end{document}